\documentclass[11pt]{amsart}
\usepackage{amsmath, amsthm, amssymb, amscd}
\usepackage{graphics}

\title[numerical radius Haagerup norm]
{The numerical radius Haagerup norm
and Hilbert space square factorizations}

\author[T. Itoh]{Takashi Itoh$^*$}
\address{$^*$Department of Mathematics, Gunma University,
         Gunma 371-8510, Japan}
\email{itoh$@$edu.gunma-u.ac.jp}
\author[M. Nagisa]{Masaru Nagisa$^{**}$}
\address{$^{**}$Department of Mathematics and Informatics, Chiba
University,
         Chiba 263-8522, Japan}
\email{nagisa$@$math.s.chiba-u.ac.jp}

\newtheorem{thm}{Theorem}[section]
\newtheorem{prop}[thm]{Proposition}

\newtheorem{cor}[thm]{Corollary}

\theoremstyle{definition}

\theoremstyle{remark}

\newtheorem{rem}[thm]{Remark}

\date{}

\begin{document}

\begin{abstract}
We study a factorization of bounded linear maps 
from an operator space $A$ to its dual space $A^*$.
It is shown that 
$T : A \longrightarrow A^*$ factors through a pair of 
a column Hilbert spaces $\mathcal{H}_c$ and its dual space
if and only if $T$ is a bounded linear form on 
$A \otimes A$ by the canonical identification
equipped with a numerical radius type Haagerup norm.
As a consequence, we characterize a bounded linear map 
from a Banach space to its dual space,
which factors through a pair of Hilbert spaces.
\end{abstract}

\maketitle




\section{Introduction}\label{intro}

The factorization through a Hilbert space of a linear map 
plays one of the central roles in the Banach space theory 
(c.f. \cite{pisier1}).
Also in the $C^*$-algebra and the operator space theory,
many important factorization theorems have been proved 
related to the Grothendieck type inequality
in several situations  \cite{haagerup},  \cite{effrosruan1},
\cite{pisieroh},  \cite{pisiershlyakhtenko}.

Let $\alpha$ be a bounded linear map from $\ell^1$ to $\ell^\infty$,
$\{e_i\}_{i=1}^\infty$ the canonical basis of $\ell^1$,
and $\mathbb{B}(\ell^2)$ the bounded operators on $\ell^2$.  
We regard $\alpha$ as the
infinite dimensional matrix $[\alpha_{ij}]$ where 
$\alpha_{ij}=\langle e_i, \alpha(e_j) \rangle$.
The Schur multiplier $S_\alpha$ on $\mathbb{B}(\ell^2)$ is defined by
$S_\alpha(x)=\alpha \circ x$ for $x=[x_{ij}] \in \mathbb{B}(\ell^2)$ 
where $\alpha \circ x$ is the Schur
product $[\alpha_{ij} x_{ij}] $. Let $w(\cdot)$ be the numerical
radius norm on $\mathbb{B}(\ell^2)$.
In \cite{itohnagisa2}, it was shown that 
$$\| S_\alpha \|_w = \sup_{x \neq 0} \frac{w(\alpha \circ x)}{w(x)} \le 1$$ 
if and only if 
$\alpha$ has the following factorization with $\|a\|^2 \|b\| \le 1$:

\[  \begin{CD}  \ell^1  @>{\alpha}>>  \ell^\infty   \\
                @V{a}VV            @AA{a^t}A     \\
                \ell^2  @>>{b}>    {\ell^2}^*
    \end{CD}
\]
where $a^t$ is the transposed map of $a$.

Motivated by the above result,  we will show a square factorization
theorem of a bounded linear map through a pair of column Hilbert 
spaces $\mathcal{H}_c$ between an operator space 
and its dual space.
More precisely, let us suppose that $A$ is an operator space in
$\mathbb{B}(\mathcal{H})$ and  $A \otimes A$ is the algebraic 
tensor product.
We define the numerical radius Haagerup norm of an element 
$u \in A \otimes A$ by
\[  \| u \|_{wh} =  \inf \{ \frac{1}{2} \|[x_1, \dots,  x_n,  
                                        y_1^*, \dots, y_n^*]\|^2 
    \mid  u = \sum_{i=1}^n x_i \otimes y_i  \} .
\]
Let $T : A \longrightarrow A^*$ be a bounded linear map.
We show that $T : A \longrightarrow A^*$ has an extention $T'$ which 
factors through a pair of column Hilbert spaces 
$\mathcal{H}_c$ so that
\[  \begin{CD}  C^*(A)  @>{T'}>>  C^*(A)^*   \\
                @V{a}VV            @AA{a^*}A     \\
                \mathcal{H}_c  @>>{b}>    \mathcal{H}_c
    \end{CD} 
\]
with  $\inf \{ \| a \|_{cb}^2 \| b \|_{cb} \mid T'= a^* b a \} \le 1$
if and only if 
$T \in (A \otimes_{wh} A)^*$ with $\| T \|_{{wh}^*}  \le 1$
by the natural identification
$\langle x, T(y) \rangle = T(x \otimes y)$ for $x, y \in A$.

We also study a variant of the numerical radius Haagerup norm 
in order to get the factorization without using the $*$ structure.

As a consequence, the above result and/or the variant read a square
factorization of a bounded linear map through a pair of Hilbert spaces 
from a Banach space $X$ to its dual space $X^*$.
The norm on $X \otimes X $ corresponding to the numerical radius
Haagerup norm is as follows:
\[  \| u \|_{wH} =  \inf\{ \sup \{ 
     (\sum_{i=1}^n|f(x_i)|^2)^{\frac{1}{2}}
     (\sum_{i=1}^n |f(y_i)|^2)^{\frac{1}{2}} \} \} ,  
\]
where the supremum is taken over all $f \in X^*$ with $\|f\| \le 1$ 
and the infimum is taken over all representation
$ u = \sum_{i=1}^n x_i \otimes y_i \in X \otimes X$.

The norm $\| \ \|_{wH}$ is equivalent to the norm $\| \ \|_{H}$ 
(see Remark 4.4), which was introduced by Grothendieck 
in \cite{grothendieck}.
However $\| \ \|_{wH}$ will give another point of view to 
a bounded linear map when we consider the factorization
through Hilbert spaces.
Let $\pi_2(a)$ be the 2-summing norm (c.f.\cite{pisier1} or see section 4) 
of a linear map $a$ from $X$ to $\mathcal{H}$.
We show that $T : X \longrightarrow X^*$ has the factorization
\[  \begin{CD}  X  @>{T}>>  {X^*}   \\
                @V{a}VV            @AA{a^t}A     \\
                \mathcal{H}  @>>{b}>    \mathcal{{H}^*}
    \end{CD} 
\]
with $\inf  \{\pi_2(a)^2 \| b \| \mid T= a^t b a \} \le 1$
if and only if $T \in (X \otimes_{wH} X)^*$ with 
$\| T \|_{{wH}^*}  \le 1$.
Moreover we characterize a linear map $X \longrightarrow X^*$ 
which has a square factorization by 
a Lindenstrauss and Pelczynski type condition 
(c.f. \cite{lp} or see Remark 4.4)

We refere to \cite{effrosruan2}, \cite{paulsen}, \cite{pisier3}  for
background on operator spaces, 
\cite{pisier1}, \cite{pisier2}  for factorization through 
a Hilbert space, and  \cite{paulsensuen}, \cite{suen1}, \cite{suen}, \cite{suen3}
for completely bounded maps related to the numerical radius norm.

\vspace{20pt}


\section{Factorization on operator spaces}

Let $\mathbb{B}(\mathcal{H})$ be the space of all bounded operators 
on a Hilbert space $\mathcal{H}$. 
Throughout this paper,  
let us suppose that $A$ and $B$ are operator spaces
in $\mathbb{B}(\mathcal{H})$.
We denote by $C^*(A)$ the $C^*$-algebra in $\mathbb{B}(\mathcal{H})$
generated by the operator space $A$.
We define the numerical radius Haagerup norm of an element 
$u \in A \otimes B$  by
\[  \| u \|_{wh} =  \inf \{\frac{1}{2}
                \|[x_1, \dots, x_n, y_1^*, \dots, y_n^*]\|^2
       \mid  u = \sum_{i=1}^n x_i \otimes y_i \} ,
\]
where $[x_1, \dots, x_n, y_1^*, \dots, y_n^*]\in M_{1, 2n}(C^*(A+B))$,
and denote by $A \otimes_{wh} B$ the completion of $A \otimes B$ 
with the norm $\| \ \|_{wh}$.

Recall that the Haagerup norm on $A \otimes B$ is
\[  \| u \|_{h} =  \inf \{ \|[x_1, \dots , x_n]\| 
                           \| [y_1, \dots ,y_n]^t \|
        \mid   u = \sum_{i=1}^n x_i \otimes y_i \},
\]
where $[x_1, \dots, x_n] \in M_{1,n}(A)$ and 
$[y_1, \dots, y_n]^t \in M_{n,1}(B)$.

By the identity
$$   \inf_{\lambda > 0} \frac{\lambda \alpha + \lambda^{-1} \beta}{2} 
    = \sqrt{\alpha \beta}  \eqno{(\star)}
$$
for positive real numbers $\alpha, \beta \ge 0$, 
the Haagerup norm can be rewritten as
\[  \| u \|_{h} =  \inf \{\frac{1}{2}
           (\|[x_1, \dots , x_n]\|^2 + \|[y_1^*, \dots, y_n^*] \|^2) 
    \mid   u = \sum_{i=1}^n x_i \otimes y_i \}.
\]

Then it is easy to check that
\[  \frac{1}{2} \| u \|_h \le \| u \|_{wh} \le \| u \|_h  \]
and $\| u \|_{wh}$ is a norm.
We use the notation  $x \alpha \odot y^t$ for
$\sum_{i=1}^n \sum_{j=1}^m x_i \alpha_{ij} \otimes y_j $,
where $x= [x_1, \dots, x_n] \in M_{1,n}(A)$, 
$\alpha =[\alpha_{ij}] \in M_{n,m}(\mathbb{C})$
and $ y^t =[y_1, \dots, y_m]^t \in M_{m,1}(B)$.  
We note the identity
$x \alpha \odot y^t = x \odot \alpha y^t$.

First we show that the numerical radius Haagerup norm has the
injectivity.

\vspace{20pt}

\begin{prop}
Let $A_1 \subset A_2$ and $B_1 \subset B_2$ be operator spaces in
$\mathbb{B}(\mathcal{H})$.
Then the canonical inclusion $\Phi$ of $A_1 \otimes_{wh} B_1$ 
into $A_2 \otimes_{wh} B_2$ is isometric.
\end{prop}

\begin{proof}
The inequality $\| \Phi(u) \|_{wh} \le \| u \|_{wh} $ is trivial.  
To get the converse inequality,
let $u = \sum_{i=1}^n x_i \otimes y_i \in A_1 \otimes B_1$.
We may assume that
$\{y_1, \dots, y_k \} \subset B_2 $
is linearly  independent and  there exists
an $n \times k$ matrix of scalars 
$L \in M_{nk}(\mathbb{C})$ such that
$[y_1, \dots, y_n]^t = L[y_1, \dots , y_k]^t$.
We put $z^t=[y_1, \ldots y_k]^t$.
Then we have
\begin{align*}
  u & = x \odot y^t = x \odot Lz^t  \\
    & = xL(L^*L)^{-1/2}\odot (L^*L)^{1/2}z^t
\end{align*}
and
\[  \| [ xL(L^*L)^{-1/2}, ((L^*L)^{1/2}z^t)^* ]\| \le
    \| [x, (y^t)^*] \|.  
\]
So we can get a representation
$u=[x'_1, \ldots , x'_k]\odot [y'_1, \ldots, y'_k]^t$ with
\[  \| [x'_1, \ldots , x'_k, {y'_1}^*, \ldots, {y'_k}^*] \| \le
    \| [x, (y^t)^*] \| 
\]
and $ \{ y'_1, \ldots ,y'_k\}$ is linearly independent.
This implies that $x'_1, \ldots , x'_k \in A_1$.

Applying the same argument for $\{ x'_1, \ldots, x'_k\}$ instead of
$\{y_1, \ldots , y_n \}$,  we can get a representation
$u=[x"_1, \ldots , x"_l]\odot [y"_1, \ldots, y"_l]^t$ with
\[  \| [x"_1, \ldots , x"_l, y"_1 ^*, \ldots, y"_l ^*] \| \le
    \| [x, (y^t)^*] \| 
\]
and $ x"_i \in A_1$ and $y"_i \in B_1$.
It follows that  $\| \Phi(u)\|_{wh} \ge \| u\|_{wh}$.
\end{proof}

\vspace{20pt}

We also define a norm of an element $u \in C^*(A) \otimes C^*(A)$ by
\[  \| u \|_{Wh} =  \inf \{\|[x_1, \dots , x_n]^t \|^2w(\alpha) 
    \mid   u = \sum x_i^* \alpha_{ij} \otimes x_j \},
\]
where $w(\alpha)$ is the numerical radius norm of 
$\alpha =[\alpha_{ij}]$ in $M_n(\mathbb{C})$.

$A \otimes_{Wh} A$ is defined as the closure of 
$A \otimes A$ in $C^*(A) \otimes_{Wh} C^*(A)$.

\vspace{20pt}

\begin{thm}
Let $A$ be an operator space in $\mathbb{B}(\mathcal{H})$.
Then $A \otimes_{wh} A = A \otimes_{Wh} A$.
\end{thm}

\begin{proof}
By Proposition 2.1 and the definition of $A \otimes_{Wh} A$, 
it is sufficient to show that
$C^*(A) \otimes_{wh} C^*(A) = C^*(A) \otimes_{Wh} C^*(A)$.

Given $u = \sum_{i=1}^n x_i \otimes y_i \in C^*(A) \otimes C^*(A)$, 
we have
\[ u= [x_1,\ldots ,x_n, y_1^*,\ldots ,y_n^*]
   \begin{bmatrix} 0_n & 1_n \\ 0_n & 0_n \end{bmatrix} \odot
   [x_1^*,\ldots ,x_n^*, y_1,\ldots ,y_n]^t . \]
Since 
$w( \begin{bmatrix} 0_n & 1_n \\ 0_n & 0_n \end{bmatrix})=\frac{1}{2}$,
then $ \| u \|_{wh} \ge \| u \|_{Wh}$.

To establish the reverse inequality,
suppose that
$u = \sum_{i,j=1}^n x_i^* \alpha_{ij} \otimes x_j 
     \in C^*(A) \otimes C^*(A)$ 
with $w(\alpha) =1$ and $\|[x_1, \dots, x_n]^t\|^2 =1$.
It is enough to see that there exist 
$c_i, d_i \in C^*(A) (i=1, \dots, m)$
such that $u = \sum_{i=1}^m c_i \otimes d_i$  
with $\|[c_1, \dots, c_m, d_1^*, \dots, d_m^*] \|^2  \le 2$.
By the assumption $w(\alpha) =1$ and Ando's Theorem \cite{ando}, 
we can find a self-adjoint matrix $\beta \in M_n(\mathbb{C})$
for which 
$ P = \begin{bmatrix} 1+\beta & \alpha \\
                 \alpha^* & 1-\beta \end{bmatrix}$
is positive definite in $M_{2n}(\mathbb{C})$. 

Set $[c_1, \dots, c_{2n}] =[x_1^*, \dots, x_n^*, 0, \dots,0]
      P^{\frac{1}{2}}$
and
$[d_1, \dots d_{2n}]^t = P^{\frac{1}{2}} 
    [0, \dots, 0,$ $x_1, \dots, x_n]^t$. 
We note that $u = [c_1, \dots, c_{2n}] \odot [d_1, \dots, d_{2n}]^t $.
Then we have
\begin{align*}
    & \|[[c_1, \dots, c_{2n}, d_1^*, \dots, d_{2n}^*] \|^2 \\
  = & \| [x_1^*, \dots, x_n^*, 0, \dots, 0]P
              [x_1, \dots , x_n , 0, \dots , 0]^t \\
    & \qquad \qquad + [0, \dots, 0, x_1^*, \dots, x_n^*]P
                    [0, \dots, 0,x_1, \dots, x_n]^t  \|  \\
  = &  \| [x_1^*, \dots, x_n^*](1+\beta) [x_1, \dots, x_n]^t \\
    &  \qquad \qquad +
          [x_1^*, \dots, x_n^*](1-\beta)[x_1, \dots, x_n]^t \| \\
  = &  2\|[x_1, \dots, x_n]^t \|^2  = 2.
\end{align*}

\end{proof}

\vspace{20pt}

We recall the column (resp. row) Hilbert space $\mathcal{H}_c $(resp.
$\mathcal{H}_r$) for a
Hilbert space $\mathcal{H}$.
If $\xi = [\xi_{ij}] \in M_n(\mathcal{H})$,  then we define a map
$C_n(\xi)$ by
$$
C_n(\xi) : \mathbb{C}^n \ni [\lambda_1, \dots, \lambda_n] \longmapsto
[\sum_{j=1}^n \lambda_j \xi_{ij}]_i \in \mathcal{H}^n
$$
and denote the column matrix norm by $\|\xi \|_c = \|C_n(\xi)\|$. This
operator
space
structure on $\mathcal{H}$ is called
the column Hilbert space and denoted by $\mathcal{H}_c$.

To consider the row Hilbert space, let  $ {\overline{\mathcal{H}}}$ be
the conjugate Hilbert space for $\mathcal{H}$.  We define a map
$R_n(\xi)$
by
$$
R_n(\xi) : {\overline{\mathcal{H}}}^n \ni [{\overline{\eta}}_1, \dots,
{\overline{\eta}}_n] \longmapsto
[\sum_{j=1}^n (\xi_{ij} | \eta_j) ]_i  \in {\mathbb{C}}^n
$$
and  the row matrix norm by $\|\xi \|_r = \|R_n(\xi)\|$.
This operator space structure on $\mathcal{H}$ is called
the row Hilbert space and denoted by $\mathcal{H}_r$.

Let  $a : C^*(A) \longrightarrow \mathcal{H}_c$ be a completely bounded
map.
We define a map $d : C^*(A) \longrightarrow \overline{\mathcal{H}}$
by $d(x) = \overline{ a(x^*)}$ .  It is not hard to
check that $d : C^*(A) \longrightarrow \overline{\mathcal{H}}_r$
is completely bounded and $\|a\|_{cb} = \|d\|_{cb}$ when we introduce
the
row
Hilbert space structure
to $\overline{\mathcal{H}}$.
In this paper, we define the adjoint map $a^*$ of $a$
by the transposed map of $d$, that is,
$d^t :(( \overline{\mathcal{H}})_r)^*
= ((\mathcal{H}^*)_r)^* = (\mathcal{H}^{**})_c =\mathcal{H}_c
\longrightarrow
C^*(A)^*$
(c.f. \cite{effrosruan1}).
More precisely,
we define
$$
\langle a^*(\eta ), x \rangle = \langle \eta, d(x) \rangle = (\eta |
a(x^*))
\qquad {\text{for}} \qquad \eta \in \mathcal{H}, x \in C^*(A).
$$

A linear map $T :  A \longrightarrow A^*$ can be identified with
the bilinear form $ A \times A \ni (x, y) \longmapsto \langle x,
T(y)\rangle \in
\mathbb{C}$
and also the linear form $ A \otimes A \longrightarrow \mathbb{C}$.
We use $T$  also to denote both of the bilinear form and the linear
form,
and $\| T \|_{{\beta}^*}$ to denote the
norm when $A \otimes A$ is equipped with a norm $\| \ \|_\beta$.

We are going to prove the main theorem. The proof will be given
by the similar way as in the case of the
original Haagerup norm in \cite{effroskishimoto}.

\vspace{20pt}

\begin{thm}
Suppose that $A$ is an operator space in $\mathbb{B}(\mathcal{H})$, 
and that  $T : A \times A \longrightarrow \mathbb{C}$ is bilinear.
Then the following are equivalent:
\begin{enumerate}
  \item[(1)] $\|T\|_{wh^*}  \le 1$.
  \item[(2)] There exists a state $p_0$ on $C^*(A)$  such 
   \[   |T(x, y)| \le p_0(xx^*)^{\frac{1}{2}}  p_0(y^*y)^{\frac{1}{2}}
        \qquad \text{ for } x, y \in A.
   \]
  \item[(3)] There exist a $*$-representation 
      $\pi : C^*(A) \longrightarrow \mathbb{B}(\mathcal{K})$,
      a unit vector $\xi \in \mathcal{K}$ and 
      a contraction $b \in \mathbb{B}(\mathcal{K})$ such that
    \[  T(x, y) = ( \pi(x) b \pi(y) \xi \mid \xi)  
        \qquad \text{ for } x, y \in A.
    \]
  \item[(4)] There exist an extension 
      $T' : C^*(A) \longrightarrow C^*(A)^* $ of $T$ 
      and completely bounded maps 
      $a: C^*(A) \longrightarrow \mathcal{K}_c$,
      $b: \mathcal{K}_c \longrightarrow \mathcal{K}_c$ such that
    \begin{gather*}  \begin{CD}  C^*(A)  @>{T'}>>  C^*(A)^*   \\
                                 @V{a}VV            @AA{a^*}A     \\
                                 \mathcal{K}_c  @>>{b}>    \mathcal{K}_c
                     \end{CD} \\
        \text{i.e., } T' = a^* b a \; \text{with} \; 
        \|a\|_{cb}^2\|b\|_{cb}\le 1.
    \end{gather*}
\end{enumerate}

\end{thm}

\vspace{20pt}

\begin{proof}
(1)$\Rightarrow$(2)
By Proposition 2.1, we can extend $T$ 
on $C^*(A) \otimes_{wh} C^*(A)$  and also denote it by $T$.
We may assume $\|T \|_{wh^*} \le 1$.  
By the identity $(\star)$, 
it is sufficient to show the existence of a state 
$p_0 \in S(C^*(A))$ such that
\[  |T(x, y)| \le \frac{1}{2} p_0(xx^* + y^*y) \quad \text{ for }
     x, y \in C^*(A) .
\]
Moreover it is enough to find $p_0 \in S(C^*(A))$ such that
\[  {\rm Re} T(x, y) \le \frac{1}{2} p_0(xx^* + y^*y) \quad \text{ for }
    x, y \in C^*(A).
\]
Define a real valued function  
$T_{\{x_1, \dots, x_n, y_1, \dots, y_n\}}(\ \cdot \ )$  
on $S(C^*(A))$ by
\[  T_{\{x_1, \dots, x_n, y_1, \dots, y_n \}} (p) =
    \sum_{i=1}^n \frac{1}{2}p(x_i x_i^* + y_i^*y_i) 
    - {\rm Re} T(x_i, y_i), 
\]
for $x_i, y_i \in C^*(A)$.
Set
\[  \bigtriangleup = \{ T_{\{x_1, \dots, x_n, y_1, \dots, y_n \}} 
    \mid  x_i, y_i \in C^*(A), \;  n \in \mathbb{N} \}.
\]
It is easy to see that  $\bigtriangleup$ is a cone 
in the set of all real functions on $S(C^*(A))$.
Let $\bigtriangledown$ be the open cone of all strictly negative
functions on $S(C^*(A))$.
For any $x_1, \dots, x_n, y_1, \dots, y_n \in  C^*(A)$, 
there exists $p_1 \in S(C^*(A))$ such that
$p_1(\sum x_ix_i^* + y_i^*y_i) = \| \sum x_ix_i^* + y_i^*y_i \|$.
Since
\begin{align*}
    T_{\{x_1, \dots, x_n, y_1, \dots, y_n \}} (p_1)
  & = \frac{1}{2}p_1( \sum x_i x_i^* + y_i^*y_i) 
         - {\rm Re} \sum T(x_i, y_i)  \\
  & = \frac{1}{2}\| \sum x_i x_i^* + y_i^*y_i \| 
         - {\rm Re} \sum T(x_i, y_i)  \\
  & \ge \frac{1}{2}\| \sum x_i x_i^* + y_i^*y_i \| - |\sum T(x_i, y_i)| \\
  & \ge 0,
\end{align*}
then $\bigtriangleup \cap \bigtriangledown = \phi$.

By the Hahn-Banach Theorem,  there exists a measure $\mu $ on
$S(C^*(A))$ such that $\mu(\bigtriangleup) \ge 0$ 
and $\mu (\bigtriangledown) < 0$. 
So we may assume that $\mu$ is a probability measure.
Now put $p_0 = \int p d\mu(p)$.  
Since $T_{\{x, y\}} \in \bigtriangleup$, then
\[  \frac{1}{2} p_0(xx^* + y^*y) - {\rm Re} T(x, y) 
   =  \int T_{\{x, y\}}(p) d\mu(p) \ge 0.
\]

(2)$\Rightarrow$(1)  Since
\begin{align*}
   | \sum T(x_i, y_i) | & \le  \sum p_0(x_i x_i^*)^{\frac{1}{2}} 
              p_0(y_i^* y_i)^{\frac{1}{2}} \\
  & \le \frac{1}{2} \sum p_0(x_i x_i^* + y_i^*y_i ) \\
  & \le \frac{1}{2} \|[x_1, \dots, x_n, y_1^*, \dots, y_n^*] \|^2 
\end{align*}
for $x, y \in A$,
then we have that $T \in (A \otimes_{wh} A)^* $ 
with $\| T \|_{wh^*} \le 1$.

(1)$\Rightarrow$(3) 
As in the proof of the implication (1)$\Rightarrow$(2), 
we can find a state $p \in S(C^*(A))$ such that 
$|T(x, y)| \le p(x x^*)^{\frac{1}{2}} p(y^*y)^{\frac{1}{2}}$ 
for $x, y \in C^*(A)$.
By the GNS construction, we let 
$\pi : C^*(A) \longrightarrow \mathbb{B}(\mathcal{K})$
is the cyclic representation with the cyclic vector $\xi$ and
$p(x) = ( \pi (x) \xi \mid \xi)$ for $x \in C^*(A)$.
Define a sesquilinear form on $\mathcal{K} \times \mathcal{K}$ 
by $\langle \pi (y)\xi, \pi (x) \xi \rangle = T(x^*, y)$.
This is well-defined and bounded since
\[  |\langle \pi(y)\xi, \pi(x)\xi \rangle | 
    \le  p(x^* x)^{\frac{1}{2}}p(y^* y)^{\frac{1}{2}}
    =  \| \pi(x) \xi\| \| \pi(y) \xi \|.
\]
Thus there exists a contraction $b \in \mathbb{B}(\mathcal{K})$
such that $T(x^*, y) = (b \pi (y) \xi | \pi(x) \xi)$.

(3)$\Rightarrow$(4)
Set $a(x)=\pi (x) \xi$ for $x \in C^*(A)$ and consider 
the column Hilbert structure for $\mathcal{K}$.
Then it is easy to see that $a : C^*(A) \longrightarrow \mathcal{K}_c$
is a completely contraction.
Define that $T' = a^*ba$, then it turns out $T'$ is 
an extension of $T$ and $\|a\|_{cb}^2 \|b\|_{cb} \le 1$.

(4)$\Rightarrow$(1)
Since $T'(x, y) = (b a(y) | a(x^*))$ for $x, y \in C^*(A)$, 
then we have
\begin{align*}
   | \sum_{i,j =1}^n T' & (x_i^* \alpha_{ij},  x_j) |
      = | \sum_{i,j=1}^n (b \alpha_{ij} a(x_j) \mid a(x_i)) | \\
  & = | \left( \begin{bmatrix}  b & & 0 \\
                                  & \ddots & \\
                                0 &  & b \end{bmatrix}
               \begin{bmatrix}    & &  \\
                                  & \alpha_{ij} & \\
                                  &  &   \end{bmatrix}
               \begin{bmatrix}  a(x_1)\\ \vdots \\ a(x_n)
                                  \end{bmatrix}
        \mid   \begin{bmatrix}  a(x_1)\\ \vdots \\ a(x_n)
                                  \end{bmatrix}
        \right) |   \\
  & \le w \left( \begin{bmatrix}  b & & 0 \\
                                    & \ddots & \\
                                  0 &  & b \end{bmatrix}
                 \begin{bmatrix}    & &  \\
                                  & \alpha_{ij} & \\
                                  &  &   \end{bmatrix}
         \right)
       \left\| \begin{bmatrix}  a(x_1)\\ \vdots \\ a(x_n)
                                  \end{bmatrix}
       \right\|^2 \\
  & \le \|b \|_{cb}  w(\alpha)  \|a \|_{cb}^2
       \left\| \begin{bmatrix}  x_1 \\ \vdots \\ x_n
                                  \end{bmatrix}
       \right\|^2
\end{align*}
for 
$\sum_{i,,j=1}^n x_i^* \alpha_{ij} \otimes x_j 
   \in C^*(A) \otimes C^*(A)$.
At the last inequality,
we use two facts which  $w(cd) \le \|c \| w(d)$ for double commuting
operators $c, d$,  and
$\mathbb{B}(\mathcal{K}, \mathcal{K})$ is completely isometric onto
$CB(\mathcal{K}_c, \mathcal{K}_c)$.
Hence we obtain that $\| T \|_{{wh}^*} \le \|T' \|_{{wh}^*} \le 1$.
\end{proof}

\vspace{20pt}

\begin{rem}
(i) If we replace the linear map  $\langle T(x), y \rangle = T(x, y)$
with $\langle x, T(y)\rangle = T(x, y)$ in Theorem 2.3, 
then we have a factorization of $T$ through a pair of 
the row Hilbert spaces $\mathcal{H}_r$.  
More precisely, the following condition
$(4)'$ is equivalent to the conditions in Theorem 2.3.
\begin{enumerate}
  \item[$(4)'$] There exist an extension 
     $T' : C^*(A) \longrightarrow C^*(A)^*$  of $T$ and
     completely bounded maps $a: C^*(A) \longrightarrow \mathcal{K}_r$,
     $b: \mathcal{K}_r \longrightarrow \mathcal{K}_r$ such that
   \begin{gather*}
      \begin{CD}  C^*(A)  @>{T'}>>  C^*(A)^*   \\
                  @V{a}VV            @AA{a^*}A     \\
                  \mathcal{K}_r  @>>{b}>    \mathcal{K}_r
      \end{CD} \\
      \text{i.e.}, \quad T' = a^* b a \quad \text{ with } \quad
        \|a\|_{cb}^2\|b\|_{cb}\le 1. 
   \end{gather*}
\end{enumerate}

(ii) Let $\ell^2_n$ be an $n$-dimensional Hilbert space with
the canonical basis  $\{e_1, \dots, e_n \}$.
Given $\alpha : \ell^2_n \longrightarrow {\ell^2_n}$
with $\alpha(e_j) = \sum_i \alpha_{ij}e_i$,
we set the map $\dot{\alpha} : \ell^2_n \longrightarrow {\ell^2_n}^*$
by $\dot{\alpha}(e_j)  = \sum_i \alpha_{ij}{\bar{e}}_i$ 
where $\{ {\bar{e}_i} \}$ is the dual basis.
For notational convenience, we shall also denote ${\dot{\alpha}}$ by
$\alpha$.
For $\sum_{i=1}^n x_i \otimes e_i \in C^*(A) \otimes \ell_n^2$, 
we define a norm by 
$\| \sum_{i=1}^n x_i \otimes e_i \| = \| [x_1, \dots, x_n]^t \|$.
Let $T: C^*(A) \longrightarrow C^*(A)^*$ be a bounded linear map.
Consider
$  T \otimes \alpha : C^*(A) \otimes \ell^2_n \longrightarrow  
      C^*(A)^* \otimes {\ell^2_n}^*$ 
with a numerical radius type norm ${w}(\cdot)$ given by
\[   w (T \otimes \alpha) = \sup \{ |\langle \sum x_i^* \otimes e_i, 
        T \otimes \alpha (\sum x_i \otimes e_i) \rangle | \mid
         \|\sum x_i \otimes e_i \|  \le 1 \}.
\]
Then we have
\[  \sup \{ \frac{w(T \otimes \alpha)}{w(\alpha)} \mid 
         \alpha : \ell^2_n \longrightarrow \ell^2_n, \;
         n \in \mathbb{N} \}
   = \| T \|_{{wh}^*},
\]
since
$  T( \sum x_i^* \alpha_{ij} \otimes x_j) = 
      \langle \sum x_i^* \otimes e_i, 
              T \otimes \alpha ( \sum x_i \otimes e_i) \rangle $.

(iii)  Let $u = \sum x_i \otimes y_i \in C^*(A) \otimes C^*(A)$.  
It is straightfoward from Theorem 2.3 that
\[  \|u\|_{wh} = \sup w(\sum \varphi (x_i) b \varphi (y_i))  \]
where the supremum is taken over all  $*$- preserving completely
contractions $\varphi$ and  contractions $b$.
\end{rem}

\vspace{20pt}

\section{A variant of the numerical radius Haagerup norm}

In  this section, we study a factorizaion of 
$T : A \longrightarrow A^*$ through a column Hilbert space 
$\mathcal{K}_c$ and its dual operator space ${\mathcal{K}_c}^*$.
Since the arguments and proofs of this section are almost the same as
those given in section 2,
we only indicate the places where the changes are needed.

We define a variant of the numerical radius Haagerup norm of 
an element $u \in A \otimes B$  by
\[  \| u \|_{wh'} =  \inf \{\frac{1}{2}
        \|[x_1, \dots, x_n, y_1, \dots, y_n]^t \|^2 
    \mid   u = \sum_{i=1}^n x_i \otimes y_i \},
\]
where $[x_1, \dots, x_n, y_1, \dots, y_n]^t \in M_{2n, 1}(A+B)$,
and denote by $A \otimes_{wh'} B$ 
the completion of $A \otimes B$ with the norm $\| \ \|_{wh'}$.

We remark that $\| \ \|_{wh}$ and $\| \ \|_{wh'}$ are not equivalent,
since $\| \ \|_{'h}$ in \cite{itoh} is equivalent to $\| \ \|_{wh'}$ 
and $\| \ \|_h$ and $\| \ \|_{'h}$ are not equivalent \cite{itoh},
\cite{sinclair}.

\vspace{20pt}

\begin{prop}
Let $A_1 \subset A_2$ and $B_1 \subset B_2$ be operator spaces in
$\mathbb{B}(\mathcal{H})$.
Then the canonical inclusion $\Phi$ of $A_1 \otimes_{wh'} B_1$ 
into $A_2 \otimes_{wh'} B_2$ is isometric.
\end{prop}
\begin{proof}
The proof is almost the same as that given in Proposition 2.1.
\end{proof}

\vspace{20pt}

In the next theorem, we use the transposed map 
$a^t : (\mathcal{K}_c)^* \longrightarrow C^*(A)^*$ of
$a : C^*(A)^* \longrightarrow \mathcal{K}_c$ instead of 
$a^* : \mathcal{K}_c \longrightarrow C^*(A)^*$.
We note that $(\mathcal{K}_c)^* = (\overline{\mathcal{K}})_r$ 
and the relation $a$  and $a^t$ is given by
\[  \langle a^t(\bar{\eta}), x \rangle  
     = \langle \bar{\eta}, a(x) \rangle 
     = (\bar{\eta} | \overline{a(x)})_{\overline{\mathcal{K}}} \quad
     \text{ for }
     \bar{\eta} \in  \overline{\mathcal{K}}, x \in C^*(A).
\]
It seems that the fourth condition in the next theorem is
simpler than the fourth one in Theorem 2.3, since we do not use
$*$-structure.

\vspace{20pt}

\begin{thm}
Suppose that $A$ is an operator space in $\mathbb{B}(\mathcal{H})$,
and that $T : A \times A \longrightarrow \mathbb{C}$ is bilinear.  
Then the following are quivalent:
\begin{enumerate}
  \item[(1)] $\|T\|_{{wh'}^*}  \le 1$.
  \item[(2)] There exists a state $p_0$ on $C^*(A)$  such that
    \[  |T(x, y)| \le p_0(x^*x)^{\frac{1}{2}}  p_0(y^*y)^{\frac{1}{2}}
        \quad \text{ for }  x, y \in A.
    \]
  \item[(3)] There exist a $*$-representation 
    $\pi : C^*(A) \longrightarrow \mathbb{B}(\mathcal{K})$,
    a unit vector $\xi \in \mathcal{K}$ and a contraction  
    $b :\mathcal{K} \longrightarrow \overline{{\mathcal{K}}}$
    such that
    \[  T(x, y) = ( b \pi(y) \xi \mid  
        \overline{\pi(x)\xi})_{\overline{{\mathcal{K}}}} \quad
        \text{ for }  x, y \in A.
    \]
  \item[(4)] There exist a completely bounded map 
    $a: A \longrightarrow \mathcal{K}_c$ and a bounded map 
    $b: \mathcal{K}_c \longrightarrow (\mathcal{K}_c)^*$
    such that
    \begin{gather*}  \begin{CD}  A  @>{T}>>  A^*   \\
                     @V{a}VV            @AA{a^t}A     \\
                     \mathcal{K}_c  @>>{b}>    (\mathcal{K}_c)^*
                     \end{CD} \\
    \text{i.e.}, \quad T = a^t b a \quad \text{ with } \quad  
    \|a\|_{cb}^2\|b\| \le 1. 
    \end{gather*}
\end{enumerate}
\end{thm}

\begin{proof}
(1)$\Rightarrow$(2)$\Rightarrow$(3)
We can prove these implications by the similar way as in the proof of
Theorem 2.3.

(3)$\Rightarrow$(4) We note that we use the norm $\| \ \|$ for $b$ 
instead of the completely bounded norm $\| \ \|_{cb}$.

(4)$\Rightarrow$(1) For $x_i, y_i \in A$, we have
\begin{align*}
   | \sum_{i =1}^n T & (x_i,  y_i) |
      = | \sum_{i=1}^n (b a(y_i) \mid
          \overline{a(x_i)})_{\overline{\mathcal{K}}} |  \\
  & = \left| \left(
          \begin{bmatrix} 0 &           & & b&         &  \\
                            & \ddots & &  &\ddots &  \\
                            &           &0&  &          &b \\
                          0 &           & & 0&         &  \\
                            & \ddots & &  &\ddots &  \\
                            &           &0&  &          &0 
          \end{bmatrix}
          \begin{bmatrix} \overline{a(x_1)}\\  \vdots \\
                          \overline{a(x_n)}\\  a(y_1)\\
                          \vdots \\   a(y_n)
          \end{bmatrix}
          \right. \left|
          \begin{bmatrix} \overline{a(x_1)}\\  \vdots \\
                          \overline{a(x_n)}\\  a(y_1)\\
                          \vdots \\  a(y_n)
          \end{bmatrix}
       \right) \right|   \\
&\le w \left(
\left[ \begin{array}{cccccc}
0 &           & & b&         &  \\
   & \ddots & &  &\ddots &  \\
   &           &0&  &          &b \\
0 &           & & 0&         &  \\
   & \ddots & &  &\ddots &  \\
   &           &0&  &          &0 \\
\end{array} \right]
\right)
\left\|
\left[
         \begin{array}{c}
         \overline{a(x_1)}\\
           \vdots \\
          \overline{a(x_n)}\\
             a(y_1)\\
           \vdots \\
             a(y_n)\\
         \end{array}
         \right]
\right\|^2 \\
&= \frac{1}{2} \|b\| \|a\|_{cb}^2 \| [ x_1, \dots, x_n, y_1, \dots,
y_n]^t
\|^2
\\
&\le \frac{1}{2} \| [ x_1, \dots, x_n, y_1, \dots, y_n]^t \|^2 .
\end{align*}
\end{proof}

\vspace{20pt}

\section{Factorization on Banach spaces}

Let $X$ be a Banach space.  
Recall that the minimal quantization ${\rm Min}(X)$ of $X$.
Let $\Omega_X$ be the unit ball of $ X^*$, that is,
$\Omega_X = \{ f \in X^* | \ \|f \| \le 1 \}$.
For $[x_{ij}] \in M_n(X)$,  $\| [x_{ij}] \|_{{\rm min}}$ is defined by
\[  \| [x_{ij}] \|_{{\rm min}} = 
           \sup \{ \| [f(x_{ij})] \| \ | \  f \in \Omega_X \}.
\]
Then ${\rm Min}(X)$ can be regarded as a subspace 
in the $C^*$-algebra $C(\Omega_X)$ of all continuous
functions on the compact Hausdorff space $\Omega_X$.
Here we define a norm of an element $u \in X \otimes X$ by
\[  \| u \|_{wH} =  \inf\{ \sup
       \{ (\sum_{i=1}^n |f(x_i)|^2)^{\frac{1}{2}}
          (\sum_{i=1}^n |f(y_i)|^2)^{\frac{1}{2}}\}\},
\]
where the supremum is taken over all $f \in X^*$ with $\|f\| \le 1$
and the infimum is taken over all representation
$ u = \sum_{i=1}^n x_i \otimes y_i $.

\vspace{20pt}

\begin{prop}
Let $X$ be a Banach space. 
Then
\[  {\rm Min}(X) \otimes_{wh} {\rm Min}(X) 
     = {\rm Min}(X) \otimes_{wh'} {\rm Min}(X) = X \otimes_{wH} X.
\]
\end{prop}

\begin{proof}
Let $u = \sum_{i=1}^n x_i \otimes y_i \in {\rm Min}(X)$.
Then, using the identity $(\star)$, we have
\begin{align*}
 & \|u\|_{wh}
  = \inf \{\frac{1}{2} \|[ x_1, \dots, x_n, y_1^*, \dots, y_n^*] \|^2
 | \
                    u  = \sum_{i=1}^n x_i \otimes y_i \} \\
= &  \inf  \{\sup \{\frac{1}{2}
            \|[ f(x_1), \dots, f(x_n), \overline{f(y_1)}, \dots,
\overline{f(y_n)}] \|^2 \ | \ f \in \Omega_X \} \ | \
                    u  = \sum_{i=1}^n x_i \otimes y_i \} \\
= &  \inf  \{\sup  \{\frac{1}{2}(\sum_{i=1}^n |f(x_i)|^2 + |f(y_i)|^2)
 |
\ f
\in \Omega_X \} \ | \
                    u  = \sum_{i=1}^n x_i \otimes y_i \} \\
= &  \inf  \{\sup  \{(\sum_{i=1}^n |f(x_i)|^2)^{\frac{1}{2}}
(\sum_{i=1}^n
|f(y_i)|^2)^{\frac{1}{2}}  \ | \ f \in \Omega_X \} \ | \
                    u  = \sum_{i=1}^n x_i \otimes y_i \} \\
= & \|u \|_{wH}.
\end{align*}

The equality $\|u \|_{wh'}= \|u \|_{wH}$ is obtained 
by the same way as above.
\end{proof}

\vspace{20pt}

Let $T : X \longrightarrow X^*$ be a bounded linear map. 
As in Remark 2.4(ii), we consider the map
$T \otimes \alpha : 
   X \otimes \ell^2_n \longrightarrow X^* \otimes {\ell^2_n}^*$ 
and define a norm for
$\sum x_i \otimes e_i \in X \otimes \ell^2_n$
by
\[  \| \sum x_i \otimes e_i \| = 
       \sup \{ (\sum |f(x_i)|^2)^\frac{1}{2} \mid
       f \in \Omega_X \}.
\]
We note that,  given $x \in X$,  $x^*$ is regarded as
$ \langle x^*, f \rangle = \overline{f(x)}$ 
for $f \in X^*$ in the definition of
$w(T \otimes \alpha)$, that is,
\[  w(T \otimes \alpha) =
      \sup \{ | \langle \sum x_i^* \otimes e_i, T \otimes 
                \alpha( \sum x_i \otimes e_i) \rangle | \mid 
      \| \sum x_i \otimes e_i \| \le 1 \}.
\]
Let $a : X \longrightarrow Y$ be a linear map between Banach spaces.
$a$ is called a 2-summing operator if there is a constant $C$ which satisfies
the inequality
$$
(\sum\|a(x_i)\|^2)^{\frac{1}{2}} \le C
\sup \{ (\sum |f(x_i)|^2)^{\frac{1}{2}} \ | \ f \in \Omega_X \}
$$
for any finite subset $\{x_i \} \subset X$. $\pi_2(a)$ is the smallest constant
of $C$, and is called the 2-summing norm of $a$.
The following might be well known.

\vspace{20pt}

\begin{prop}
Let $X$ be a Banach space. If a is a linear map from $X$ to $\mathcal{H}$,
then the following are equivalent:
\begin{enumerate}
  \item[(1)] $\|a : \rm{Min}(X) \longrightarrow \mathcal{H}_c \|_{cb} \le 1.$
  \item[(2)] $\|a : \rm{Min}(X) \longrightarrow \mathcal{H}_r \|_{cb} \le 1.$
  \item[(3)] $\pi_2(a : X \longrightarrow H) \le 1.$
\end{enumerate}
\end{prop}

\begin{proof}
(1) $\Rightarrow$ (3)
For any $x_1, \dots , x_n \in X$, we have
\begin{align*}
\sum_{i=1}^n \|a(x_i)\|^2 & = \| [a(x_1), \cdots,  a(x_n) ]^t \|^2 \\
& = \|a\|_{cb}^2 \|[x_1, \cdots, x_n]^t \|^2_{\text{Min}} \\
& = \| \sum_{i=1}^n x_i^*x_i \|_{\text{Min}} \\
& = \sup \{ \sum_{i=1}^n |f(x_i)|^2 \ | \ f \in \Omega_X \}.
\end{align*}

(3) $\Rightarrow$ (1)
For any $[x_{ij}] \in {M}_n({\text{Min}}(X))$,  we have
\begin{align*}
\|[a(x_{ij})]\|^2_{{M}_n(\mathcal{H}_c)} 
& = \sup \{ \sum_i\| \sum_j \lambda_ja(x_{ij}) \|^2 \ | \ \sum |\lambda_j|^2 =1 \} \\
&\le \sup \{ \pi_2(a)^2 \sup \{ \sum_i |f(\sum_j \lambda_j x_{ij}) |^2 \ 
| \ f \in \Omega_X \} \ | \ \sum |\lambda_j|^2 =1 \} \\
&\le \sup \{ \|[f(x_{ij})] \|^2 \ | \ f \in \Omega_X \} \\
& \le \|[x_{ij}] \|^2_{{M}_n({\text{Min}}(X))}.
\end{align*}

(2) $\Leftrightarrow$ (3)
It follows from the same way as above.
\end{proof}

Finally we can state the following result as a corollary of Theorem 2.3
and Theorem 3.2.

\vspace{20pt}

\begin{cor}
Suppose that $X$ is a Banach space, and that 
$T : X \longrightarrow X^*$ is a bounded linear map.
Then the following are quivalent:
\begin{enumerate}
  \item[(1)] $w(T \otimes \alpha) \le w(\alpha)$ 
    for all $\alpha : \ell^2_n \longrightarrow {\ell^2_n}$
    and $n \in \mathbb{N}$.
  \item[(2)] $\|T\|_{{wH}^*} \le 1$.
  \item[(3)] $T$ factors through a Hilbert space $\mathcal{K}$ and its
    dual space $\mathcal{K}^*$ by a 2-summing operator
    $a: X \longrightarrow \mathcal{K}$ and a bounded operator
    $b: \mathcal{K}\longrightarrow {\mathcal{K}}^*$ 
    as follows:
    \begin{gather*} \begin{CD}  X  @>{T}>>  X^*   \\
                    @V{a}VV            @AA{a^t}A     \\
                    \mathcal{K}  @>>{b}>    {\mathcal{K}}^*
                    \end{CD}   \\
    \text{i.e.}, \quad T= a^t b a \quad \text{ with } \quad  
    \pi_2(a)^2\|b\| \le 1.
    \end{gather*}
  \item[(4)] $T$ has an extention 
    $T' : C(\Omega_X) \longrightarrow C(\Omega_X)^* $ 
    which factors through a pair of Hilbert spaces 
    $\mathcal{K}$ by a 2-summing operator
    $a:C(\Omega_X) \longrightarrow \mathcal{K}$  and a bounded operator
    $b: \mathcal{K}\longrightarrow \mathcal{K}$ 
    as follows:
    \begin{gather*} \begin{CD}  C(\Omega_X)  @>{T'}>>  C(\Omega_X)^*   \\
                                @V{a}VV            @AA{a^*}A     \\
                                \mathcal{K}  @>>{b}>    \mathcal{K}
                    \end{CD}  \\
    \text{i.e.}, \quad T' = a^* b a \quad \text{ with } \quad 
    \pi_2(a)^2\|b\| \le 1.
    \end{gather*}
\end{enumerate}
\end{cor}

\begin{proof}
(1) $\Rightarrow$ (2)
Suppose that
\[   |\langle \sum_{i = 1}^m z_i^* \otimes e_i, T \otimes 
          \alpha (\sum_{i =1}^m z_i \otimes e_i) \rangle | \le 1
\]
for any $\sum_{i = 1}^m z_i  \otimes e_i \in X \otimes \ell_m^2 $ 
with $  \| \sum_{i = 1}^m z_i \otimes e_i \| \le 1$
and $\alpha \in M_n(\mathbb{C})$ with $w(\alpha) \le 1$.
It is easy to see that 
$ | \sum_{i,j = 1}^m \langle z_i^*, T(z_j) \rangle \alpha_{ij}| \le 1$, 
equivalently
$ | \sum_{i, j = 1}^m \langle z_i, T(z_j) \rangle 
    \overline{\alpha_{ij}} | \le 1$.

Given $\| \sum_{i = 1}^n x_i \otimes y_i \|_{wH} < 1$, 
we may assume that
\[  \frac{1}{2} \| [x_1, \dots, x_n, y_1, \dots, y_n]^t \|^2 \le 1. \]
Set 
\[  z_i = \begin{cases} \frac{1}{\sqrt{2}} x_i  & i =1, \dots , n \\
                        \frac{1}{\sqrt{2}} y_{i-n}  & i=n+1, \dots , 2n
          \end{cases} \quad \text{ and }
    \alpha = \begin{bmatrix} 0_n & 2\cdot 1_n \\
                             0_n & 0_n
             \end{bmatrix} .
\]
It turns out  $\| \sum_{i = 1}^{2n} z_i  \otimes e_i \| \le 1$ and
$w(\alpha) =
1$.
Then we have $|T(\sum_{i = 1}^n x_i \otimes y_i )| =| \sum_{i, j =
1}^{2n}
\langle z_i, T(z_j) \rangle \alpha_{ij} | \le 1.$
Hence $\|T\|_{wH^*} \le 1$.

(2) $\Rightarrow$ (1)
Suppose that $\|T\|_{wH^*} \le 1$. 
Then $T$ has an extension
$T' \in (C(\Omega_X) \otimes_{wh} C(\Omega_X))^*$ 
with $\| T' \|_{wh^*}\le 1$.
Given $\varepsilon >0$ and $\alpha \in M_n(\mathbb{C})$, 
there exist $x_1, \dots, x_n \in C(\Omega_X)$ such that
$\| \sum_{i=1}^n x_i \otimes e_i \| \le 1$ 
(equivalently $\|[x_1, \dots ,x_n]^t \| \le 1$) and
$ w(T' \otimes \alpha ) - \varepsilon < 
    | \sum_{i, j =1}^n  \langle x_i^*, T'(x_j)  \rangle \alpha_{ij} |.$
Hence we have
\begin{align*}
w(T \otimes \alpha)
& \le w(T' \otimes \alpha) \\
& < |T'(\sum_{i,j=1}^n x_i^* \alpha_{ij} \otimes x_j)| + \varepsilon \\
& \le \| [x_1, \dots , x_n]^t \|^2 w(\alpha) + \varepsilon \\
& \le w(\alpha) + \varepsilon.
\end{align*}

(2) $\Leftrightarrow$ (3)
It is straightforward from Theorem 3.2 and Proposition 4.1, 4.2.

(2) $\Leftrightarrow$ (4)
It is straightforward from Theorem 2.3 and Proposition 4.1,4.2.
\end{proof}

\vspace{20pt}

\begin{rem}
Here we compare the above corollary with the classical factorization
theorems through a Hilbert space.
Let $X$ and $Y$ be Banach spaces.
Grothendieck introduced the norm $\| \ \|_H$  on $X \otimes Y$ in
\cite{grothendieck} by
$$
\| u \|_{H} =  \inf\{ \sup \{ (\sum_{i=1}^n |f(x_i)|^2)^
{\frac{1}{2}}(\sum_{i=1}^n |g(y_i)|^2)^{\frac{1}{2}}\}\}
$$
where the supremum is taken over all $f\in X^*, g \in Y^*$ with 
$\|f\|, \|g\| \le 1$
and the infimum is taken over all representation
$ u = \sum_{i=1}^n x_i \otimes y_i \in X \otimes Y$.
In \cite{lp}, Lindenstrauss and Pelczynski characterized the factorization by using
$T \otimes \alpha : X \otimes \ell_n^2 \longrightarrow 
   Y \otimes \ell_n^2$
for $T : X \longrightarrow Y$, however the norm on $X \otimes \ell^2$ is
slightly different from the one in this paper.
Their theorems with a modification are summarized
for a bounded linear map $T : X \longrightarrow Y^*$ as follows:

\vspace{20pt}

The following are equivalent:
\begin{enumerate}
  \item[(1)] $\|T \otimes \alpha \| \le \| \alpha \| $ 
     for all $\alpha : \ell^2_n \longrightarrow {\ell^2_n}$
     and $n \in \mathbb{N}$.
  \item[(2)] $\|T\|_{{H}^*} \le 1$.
  \item[(3)] $T$ factors through a Hilbert space $\mathcal{K}$ by a
     2-summing operator $a: X \longrightarrow \mathcal{K}$ and 
     $b: \mathcal{K}\longrightarrow {Y}^*$ whose transposed $b^t$ is 2-summing as follows:

\setlength{\unitlength}{0.7mm}
\begin{picture}(60, 40)(-40, 20)
\put(20, 17){\makebox(20,10){$\mathcal{K}$}}
\put(0,  40){\makebox(20,10){$X$}}
\put(40, 40){\makebox(20, 10){$Y^*$}}
\put(14, 39){\vector(1, -1){11}}
\put(34, 28){\vector(1, 1){11}}
\put(20, 45){\vector(1, 0){20}}
\put(29, 50){$T$}
\put(10, 30){$a$}
\put(46, 30){$b$}
\end{picture}
\[ \text{i.e.}, \quad T = b a \quad \text{ with } \quad  
   \pi_2(a)\pi_2(b^t) \le 1.  \]
\end{enumerate}
\end{rem}

\end{document}